\documentclass[12pt,oneside,reqno]{amsart}

\usepackage{amsthm, amssymb}
\usepackage{enumerate}
\usepackage[colorlinks]{hyperref}

\newtheorem{Exm}{Example}
\newtheorem{thm}{Theorem}[section]

\newtheorem{defi}[thm]{Definition}
\newtheorem{lem}[thm]{Lemma}
\newtheorem{prop}[thm]{Proposition}
\newtheorem{corl}[thm]{Corollary}

\allowdisplaybreaks
\title[Continuous Field of Orlicz space on groupoids]{Continuous Field of Orlicz space on locally compact groupoids and related results}
\author{K. N. Sridharan}
\address{K. N. Sridharan,\newline\indent Department of Mathematics,\newline\indent Indian Institute of Technology Delhi,\newline\indent New Delhi - 110016, India.}
\email{sreedu242@gmail.com}
\author{N. Shravan Kumar}
\address{N. Shravan Kumar,\newline\indent Department of Mathematics,\newline\indent Indian Institute of Technology Delhi,\newline\indent New Delhi - 110016, India.}
\email{shravankumar.nageswaran@gmail.com}

\begin{document}

\begin{abstract}

 Let $G$ be a locally compact second countable groupoid with a fixed Haar system $\lambda=\{\lambda^{u}\}_{u\in G^{0}}$ and $(\Phi,\Psi)$ be a complementary pair of $N$-functions satisfying the $\Delta_{2}$-condition. In this paper, we introduce the continuous field of Orlicz space $(L^{\Phi}_{0},\Delta_{1})$ over the unit space $G^0$ and provide a sufficient condition for the space of continuous sections vanishing at infinity, denoted $E^{\Phi}_{0}$, to form a Banach algebra under a suitable convolution. Additionally, we establish the condition for a closed $C_{b}(G^{0})$-submodule $I$ of $E^{\Phi}_{0}$ to be a left  ideal. Furthermore, we present a groupoid analogue of the characterization of the space of convolutors of Morse-Transue space for locally compact groups. 
\end{abstract}


\maketitle
\section{Introduction}
Harmonic analysis of groupoids has been extensively studied over the past two to three decades, with many classical results from group theory generalized to the groupoid setting. Some of the key contributions to the field include the construction of groupoid $C^{*}$-algebras \cite{renault2006groupoid} and Renault's Disintegration Theorem \cite{renault1987representation}. In 1997, Renault \cite{renault1997fourier} introduced and investigated the Fourier algebra on measured groupoids, and later, Paterson \cite{paterson2004fourier} extended this study to locally compact groupoids. Additionally, further results related to continuous representations of locally compact groupoids and the Fourier transform of compact groupoids can be found in \cite{amini2007tannaka,amini2010tannaka,bos2007groupoids}.
    
     The algebraic structure of the Orlicz space on locally compact groups is discussed in \cite{rao2001convolutions}, and specifically for locally compact abelian groups in \cite{hudzik1987some}. Additionally, a sufficient condition for the weighted Orlicz space on locally compact groups to form a Banach algebra is provided in \cite{osancliol2015weighted}. In 2019, Kumar et al. \cite{kumar2019orlicz} generalized these results to hypergroup $K$, presenting a sufficient condition for the Orlicz space to form a Banach algebra. This paper also proves the existence of an approximate identity and provides a necessary and sufficient condition for a closed subspace to be a left ideal of the Banach algebra $L^{\Phi}(K)$.
     
      It is well known that the group von Neumann algebra $VN(G)$ can be identified with the dual of the Fourier algebra $A(G)$ for any locally compact group $G$. In 2004, Paterson \cite{paterson2004fourier} generalized this result to the groupoid setting. Herz introduced and studied the $p$-analogue of the Fourier algebra, denoted by $A_p(G)$, for any locally compact group $G$. He proved that the dual of $A_p(G)$ can be identified with a subspace of the convolutors in $L^p(G)$, denoted as $PM_p(G)$. Furthermore, in 2013, Aghababa \cite{aghababa2013space} introduced and studied the convolutor space $C_{V_{\Phi}}(G)$ of the Morse-Transue space $M^{\Phi}(G)$ and showed that $C_{V_{\Phi}}(G)$ can be identified with the dual of a space denoted by $\check{A}_{\Phi}(G)$. This was inspired from the work of M. Cowling \cite{cowling1998predual} in $L_p(G)$ spaces, where $1 < p < \infty$. In 2019, Kumar et al. \cite{kumar2019orlicz} generalized this result to hypergroups.
     
    Motivated by these works, we introduce a continuous field of Orlicz spaces over the unit space $G^0$ of a locally compact groupoid $G$ and investigate the groupoid analogue of various properties of Orlicz spaces on locally compact groups and hypergroups. In this paper, we focus on the continuous sections vanishing at infinity of a field of Orlicz spaces over the unit space of groupoid and explore some of their properties. 
    
    The paper is organized as follows: In Section 2, we begin with some basic definitions and results related to locally compact groupoids and Orlicz spaces. In Section 3, we introduce a continuous field of Orlicz spaces over the unit space of the groupoid, denoted by $G^0$. Using the properties of the Haar system, we prove that the family ${(L^{\Phi}(G^{u}), \lambda^{u})}_{u \in G^0}$ forms a continuous field of Banach spaces $(L^{\Phi}_0, \Delta_1)$ or a continuous Banach bundle $(L_0^{\Phi}, p, G^0)$. We further shows that the continuous field of Orlicz spaces $\left(L_{0}^{\Phi}, \Delta_{1}\right)$ and $\left(L^{\Phi}, \Delta_{1}^{\prime}\right)$ with respect to both Gauge and Orlicz norm  respectively are isomorphic.

Next, we define the concept of a strongly continuous representation of the groupoid $G$ over a continuous field of Banach spaces and show that the left-regular representation of $G$ over $(L^{\Phi}_0, \Delta_1)$ is a strongly continuous representation. Additionally, we discuss the algebraic structure of the space of continuous sections vanishing at infinity, $E^{\Phi}_0$, and provide a sufficient condition for $E^{\Phi}_0$ to be a Banach algebra. Furthermore, we establish the existence of approximate identity and show that the closed $C_{b}\left(G^{0}\right)$-submodule $I$ of the Banach algebra $E_{0}^{\Phi}$ is a left ideal if and only if the subbundle $(I_{0}^{\Phi},p_{I_{0}^{\Phi}}, G^{0})$  is invariant under left-regular representation of $G$. 
    
       As previously noted, Aghababa \cite{aghababa2013space} introduced and studied the convolutor space $C_{V_{\Phi}}(G)$ of the Morse-Transue space $M^{\Phi}(G)$ and showed that $C_{V_{\Phi}}(G)$ can be identified with the dual of a space denoted  $\check{A}_{\Phi}(G)$. In Section 4, we extend this result to the groupoid context. Unlike the case of groups, for a locally compact groupoid $G$, we prove that the convolutor space $C_{V_{\Phi}}(G)$ is a closed subspace of the space of all bounded right-module maps  denoted $B_{D}^{\Phi}\left(\check{A}_{\Phi}(G), C_{0}(G^{0})\right)$. More specifically, we show that the closed subspace $\check{A}_{\Phi}(G)^{\prime}$ of the space of all bounded right-module maps from the $C_{0}(G^{0})$-module $\check{A}_{\Phi}(G)$ to $C_0(G^0)$ can be identified with the space of convolutors $C_{V_{\Phi}}(G)$. 

\section{Preliminaries} 

    Let us begin this section by recalling the notion of Groupoids.
    
    A groupoid is a set $G$ endowed with a product map $G^{2}\to G: (x,y)\to xy$, where  $G^{2}$ is a subset of $G\times G$ called the set of composable pairs, and an inverse map  $G\to G:x\to x^{-1}$ such that the following relations are satisfied: 
    \begin{enumerate}[(i)]
        \item $(x^{-1})^{-1}=x$,
        \item $(x,y),(y,z)\in G^{2}$ implies $(xy,z),(x,yz)\in G^{2}$ and $(xy)z=x(yz)$,
        \item $(x^{-1},x)\in G^{2}$ and if $(x,y)\in G^{2}$, then $x^{-1}(xy)=y$,
        \item $(x,x^{-1})\in G^{2}$ and if $(z,x)\in G^{2}$, then $(zx)x^{-1}=z.$
    \end{enumerate}

    If $x \in G,~ d,r: G\to G$, defined as $d(x) = x^{-1}x$ and $r(x)= xx^{-1}$ are its  domain and range maps respectively. The set $G^{0}$, the image of range and domain maps is called the unit space of $G$. Its elements are units because $xd(x) =x =r(x)x$. For $u \in G^{0}, G^{u} = r^{-1}(u)$ and $G_{u} = d^{-1}(u)$. It is also known that $(x,y)\in G^{2}$ if and only if $d(x) = r(y)$. Observe that $G = \sqcup_{u\in G^{0}}G^{u}$. For $A,B \subset G$,\[ AB=\Big\{ab \in G: a \in A,b\in B,d(a)=r(b)\Big\},~ A^{-1}=\Big\{a^{-1} \in G: a\in A\Big\}.\]
     \begin{Exm}
         \begin{enumerate}[(i)]
               \item Equivalence Relations: Let $R$ be an equivalence relation on a set $X$. Let $R^{2}=\{\big((x,y),(y,z)\big):(x,y),(y,z)\in R\}$. Define a partial product and inversion on $R$: $(x,y)(y,z)=(x,z),~ (x,y)^{-1}= (y,x)$. Then $R$ becomes a groupoid, with unit space $R^{0}=\{(x,x):x\in X\}$.
               \item Action Groupoid: Let $X$ be a set and let $\Gamma$ be a group acting on $X$ on the left. Define $G= \Gamma \times X$ and $G^{2}=\{\left((g,y),(h,x)\right): y=h \cdot x,~~ x,y \in X,~~  g,h\in \Gamma\}$ where $h \cdot x$ denote the action of $\Gamma$. Then $G$ becomes a groupoid under the given product and inversion: $(g,h \cdot x)(h,x)=(gh,x)$, ~$(g,x)^{-1}=(g^{-1},g \cdot x)$.The unit space $G^{0}=\{(e_{\Gamma},x):x\in X\}$, where $e_{\Gamma}$ is the identity element of $\Gamma$.
         \item Group bundles:The disjoint union of groups, known as a group bundle, forms a groupoid. For a collection ${G_{i}}_{i\in I}$ of groups, define $G :=\sqcup_{i\in I}(\{i\}\times G_{i})$.  The composable sets $G^{2}$ consists of those elements having same indices in their first coordinates. The product and the inverse are defined as $(i, x)(i, y)=(i, xy)$ and $(i,x)^{-1} = (i,x^{-1})$, respectively, where $xy$ and $x^{-1}$ are the product and inverse in the group $G_{i}$. The unit space $G^{0} = \{(i, e_{G_{i}}):i\in I\}$, where $e_{G_{i}}$ is the identity of $G_{i}$.
\end{enumerate} 
       \end{Exm}
    A topological groupoid consists of a groupoid $G$ and a topology compatible with the groupoid structure such that:
    \begin{enumerate}[(i)]
        \item $x\to x^{-1}:G\to G$ is continuous,
        \item$(x,y)\to xy: G^{2}\to G$ is continuous where $G^{2}$ has the induced topology from $G \times G$. 
    \end{enumerate}
       
    One observation is that the range and domain maps are continuous in a topological groupoid. Here we consider topological groupoids whose topology is Hausdorff, locally compact and we call such a groupoid $G$ a locally compact groupoid. The unit space $G^{0}$ is a locally compact Hausdorff space under the subspace topology.
    
    Like the notion of haar measure in locally compact groups, there is the notion of haar system on locally compact groupoids. The following is the definition of a haar system on locally compact groupoids.
    
    Let $G$ be a locally compact groupoid. A left Haar system for $G$ consists of a family  $\{\lambda^{u}: u\in G^{0}\}$ of positive radon measures on $G$ such that, 
\begin{enumerate}[(i)]
    \item the support of the measure $\lambda^{u}$ is $G^{u}$,
    \item for any $f\in C_{c}(G),u\to \lambda(f)(u):= \int f d\lambda^{u}$ is continuous and
    \item for any $x\in G$ and $f \in C_{c}(G)$,\[\int_{G^{d(x)}}f(xy)d\lambda^{d(x)}(y)=\int_{G^{r(x)}}f(y)d\lambda^{r(x)}(y).\]  
\end{enumerate}

    According to \cite[Proposition $2.4$]{renault2006groupoid}, if $G$ is a locally compact groupoid with a left haar system, then the range map $r$ is an open map. The same is the case for domain map $d$. This paper assumes the locally compact groupoid to be second countable.  For more details on groupoid, refer to \cite{renault2006groupoid,paterson2012groupoids}.
        
   Let $\{E(u)\}_{u\in G^{0}}$ be a family of Banach spaces. Every element $\prod_{u\in G^{0}}E(u)$, i.e. every function $\xi$ defined on $G^{0}$ such that $\xi(u)\in E(u)$ for each $u\in G^{0}$, is called a vector field \cite[Section  $10.1.1$]{MR0458185}.The vector fields are also called sections or cross-sections in various literature.
   
   A continuous field of Banach spaces over $G^{0}$ is a  family $\{E(u)\}_{u\in G^{0}}$ of Banach spaces, with a set $\Gamma\subset \prod_{u\in G^{0}}E(u)$ of vector fields such  that: 
\begin{enumerate}[(i)]
    \item $\Gamma$ is a complex linear subspace of $\prod_{u\in G^{0}}E(u)$,
    \item For every $u\in G^{0}$, the set $\xi(u)$ for $\xi \in \Gamma$ is dense in $E(u)$,
    \item For every $\xi \in \Gamma$, the function $u \to \|\xi(u)\|$ is continuous,
    \item Let $\xi \in \prod_{u\in G^{0}}E(u)$ be a vector field; if for every $u\in G^{0}$ and every $\epsilon >0$, there exists an $\xi'\in \Gamma$ such that $\|\xi(s) -\xi'(s)\| < \epsilon$ on a neighbourhood of $u$, then $\xi \in \Gamma$.
\end{enumerate}
For more details, one can  refer to \cite[Chapter $3$]{bos2007groupoids} and \cite[Section $10.1$]{MR0458185}.
    For a continuous field of Banach spaces, we can define a topology on $\mathcal{E}= \sqcup_{u\in G^{0}}E(u)$, generated by the sets of the form 
\[U(V,\xi,\epsilon)=\{h\in \mathcal{E}:\|h-\xi(p(h))\|<\epsilon,\xi\in \Gamma,p(h)\in V\},\]
    where $V$ is an open set in $G^{0}$, $\epsilon> 0$, and $p: \mathcal{E}\to G^{0}$ is the projection of the total  space $\mathcal{E}$ to base space $G^{0}$ such that fiber $p^{-1}(u) = E(u),u \in G^{0}$. This map is a surjective continuous open map under the above topology. We denote such a continuous field of  Banach spaces as $(\mathcal{E},\Gamma)$. Under the above topology, it is also referred to as a continuous Banach bundle $(\mathcal{E},p, G^{0})$ whose continuous sections are elements of $\Gamma$\cite[Definition $13.4$, Theorem $13.18$]{fell1988representations}. Here, $\mathcal{E}$ is called bundle space.
   
     For $\xi \in \Gamma$, we denote $\xi(u)$ as $\xi^{u}$ also, in the rest of the paper. Here, we call the sections or vector fields that are bounded but need not be continuous as bounded sections.
    
     A convex function $\Phi:\mathbb{R}\to [0,\infty]$ which satisfy the conditions: $\Phi(-x) =\Phi(x), \Phi(0) = 0$, and $\lim_{x\to \infty}\Phi(x)=+\infty$ is called a \textit{Young function}. For any \textit{Young function} $\Phi$, one can associate another convex function $\Psi: \mathbb{R}\to [0,\infty]$ having similar properties, called \textit{complementary function} to $\Phi$, defined by 
    \[\Psi(y)=\sup\{x|y|-\Phi(x):x\geq 0\},~y\in \mathbb{R}.\]
    
     The pair $(\Phi,\Psi)$ is called a \textit{complementary pair} of Young functions.
    
    A \textit{(nice) Young function} $\Phi$, termed as $N$-function, is a continuous Young function  such that $\Phi(x)= 0$ if and only if $x = 0$ and satisfies the following conditions:\[\lim_{x\to 0}\frac{\Phi(x)}{x} = 0~~\text{and}~~ \lim_{x\to \infty}\frac{\Phi(x)}{x} = \infty. \] 
    
   We say that a Young function $\Phi : \mathbb{R}\to \mathbb{R^{+}}$ is \textit{$\Delta_{2}$-regular} if there exist $k > 0$ and $x_{0} > 0$ such that  $$\Phi(2x)\leq k\Phi(x) \;\; \text{for all} \;\; x\geq x_{0}, \;\; \text{when} \;\; \sup_{u\in G^{0}}\lambda^{u}(G) <\infty,$$ and $$\Phi(2x)\leq k\Phi(x) \;\; \text{for all} \;\; x\geq 0, \;\; \text{when} \;\; \sup_{u\in G^{0}}\lambda^{u}(G) =\infty.$$ We write $\Phi \in \Delta_{2}$ if $\Phi$ is $\Delta_{2}$-regular. This definition is similar to the one defined in the case of locally compact groups \cite[Definition 1 (Pg. 22) and Remark (Pg. 46)]{MR1113700}.
    
    For $u \in G^{0}$, the \textit{Orlicz space}, $L^{\Phi}(G^{u})$, is defined as, \[L^{\Phi}(G^{u}) = \{f:G^u\to\mathbb{C}: f~ \text{is measurable and}~ \int_{G^{u}}\Phi(\alpha|f|)d\lambda^{u}<\infty,~\text{for some}~ \alpha>0\}.\] 
    The set $L^{\Phi}(G^{u})$ is a Banach space with respect to two norms:\\   \textit{Gauge Norm}:$\|f\|_{\Phi}^{0} =\inf\{k > 0: \int_{G^{u}}\Phi\Big(\frac{|f|}{k}\Big)d\lambda^{u}\leq 1\}$.\\   \textit{Orlicz Norm}:$\|f\|_{\Phi}=\sup\{\int_{G^{u}}|fg|d\lambda^{u}:\int_{G^{u}}\Psi(|g|)d\lambda^{u}\leq 1\}$ where $\Psi$ denotes the  complementary function to Young function $\Phi$.
    
    It is known that the two norms are equivalent. Infact, $\|\cdot\|_{\Phi}^{0}\leq \|\cdot\|_{\Phi}\leq 2\|\cdot\|_{\Phi}^{0}$ and $\|f\|_{\Phi}^{0}\leq 1$ if and only if $\int_{G^{u}}\Phi(|f|)d\lambda^{u}\leq 1$. Suppose $f$ is a real-valued measurable  function on $G^{u}$, and $\int_{G^{u}}fd\nu^{u}$ and $\int_{G^{u}}\Phi(f)d\nu^{u}$ exist, where $\nu^{u}$ is a probability  measure on $G^{u}$, then the following Jensen’s inequality\cite[Chapter $3$, Proposition $5$]{MR1113700} holds: \[\Phi\Big(\int_{G^{u}}fd\nu^{u}\Big)\leq\int_{G^{u}}\Phi(f)d\nu^{u}.\]
    The set of all complex measurable function $f$ such that: \[\int_{G^{u}}\Phi(a|f|)d\lambda^{u}<\infty,~\text{for all}~a>0.\] is denoted as $M^{\Phi}(G^{u})$. This space is called \textit{Morse-Transue space} and is a closed subspace of $L^{\Phi}(G^{u})$.The space of continuous functions on $G^{u}$ with compact support, denoted $C_{c}(G^{u})$, is dense in $M^{\Phi}(G^{u})$. If $\Phi\in \Delta_{2}$, then $L^{\Phi}(G^{u}) = M^{\Phi}(G^{u})$. Suppose  $(\Phi, \Psi)$ are complimentary pair of $N$-functions and  $(G^{u},\mathfrak{B}^{u},\lambda^{u})$ is a $\sigma$-finite measure space, where $\mathfrak{B}^{u}$ denotes the $\sigma$-algebra of subsets of $G^{u}$ on which $\lambda^{u}$ is defined, then by \cite[Chapter $4$, Theorem $7$]{MR1113700}, $M^{\Phi}(G^u)^* = L^{\Psi}(G^{u})$ and $M^\Psi(G^u)^* = L^{\Phi}(G^{u})$. For $g \in L^{\Phi}(G^{u})$,
    \begin{align*}
    &\|g\|_{\Phi}=\sup\left\{\left|\int_{G^{u}}gfd\lambda^{u}\right|:\|f\|_{\Psi}^{0}\leq 1,f\in M^{\Psi}(G^{u})\right\},\\&\|g\|_{\Phi}^{0}=\sup\left\{\left|\int_{G^{u}}gfd\lambda^{u}\right|:\|f\|_{\Psi}\leq 1,f\in M^{\Psi}(G^{u})\right\}.
    \end{align*}
    
    If $f\in L^{\Phi}(G^{u}), g \in  L^{\Psi}(G^{u})$, then $fg \in L^{1}(G^{u})$ and the following Holder’s inequality  holds: \[\int_{G^{u}}|f(t)g(t)|d\lambda^{u}(t)\leq \|f\|_{\Phi}^{0}\|g\|_{\Psi}.\] For more details, refer\cite{MR1113700}.
Throughout this paper, $C_{c}(G), C_{b}(G)$ and $C_{0}(G)$ denote the spaces of continuous functions on $G$ with compact support, bounded continuous functions on $G$ and continuous functions on $G$ that vanish at infinity, respectively.
    \section{Continuous Field of Orlicz Space}  

    Let $\{L^{\Phi}(G^{u})\}_{u\in G^{0}}$ be a family of Orlicz spaces over $G^{0}$ where $(\Phi,\Psi)$ is a complementary pair of $N$-functions with $\Phi \in \Delta_{2}$. We will prove that the above family forms a continuous field of Orlicz space over $G^{0}$. 

    \begin{thm}\label{Continuity}
        For $f\in C_{c}(G)$, the function $g : G^{0}\to \mathbb{R}$ defined as $g(u) = \|f^{u}\|_{\Phi}^{0}$,  where $f^{u} = f_{|_{G^{u}}}$, is  continuous over $G^{0}$ with compact support. 
    \end{thm}
    \begin{proof}
        Let $f\in C_{c}(G).$ The function $\Phi(|f|)$ is continuous with compact support. Hence the map $\lambda(\Phi(f))(u)=\int_{G^{u}}\Phi(|f|)d\lambda^{u}$ is continuous over $G^{0}$. Let $\epsilon>0$ and $u_{0}\in G^{0}$. Assume $\|f^{u_{0}}\|_{\Phi}^{0}\neq 0$ and we can see that $\int_{G^{u_{0}}}\Phi(\frac{|f^{u_{0}}|}{\|f^{u_{0}}\|_{\Phi}^{0}})d\lambda^{u_{0}}=1$.\\
        If $a,b \in I=(\|f^{u_{0}}\|_{\Phi}^{0}-\epsilon,\|f^{u_0}\|_{\Phi}^{0}+\epsilon)$, such that $0< a< \|f^{u_{0}}\|_{\Phi}^{0}<b$, then
    \[\lambda\Big(\Phi\Big(\frac{|f|}{a}\Big)\Big)(u_{0})>1~ \text{and}~\lambda\Big(\Phi\Big(\frac{|f|}{b}\Big)\Big)(u_{0})<1.\] By continuity of both these functions, there exist open neighbourhood $U$ of $u_{0}$ such that $\lambda\Big(\Phi\Big(\frac{|f|}{a}\Big)\Big)(u)>1$ and $\lambda\Big(\Phi\Big(\frac{|f|}{b}\Big)\Big)(u)<1$ for all $u\in U$. Thus,
    \[ a\leq \|f^{u}\|_{\Phi}^{0}\leq b ~\text{for all}~ u\in U,~\text{i.e.}~ \|f^{u}\|_{\Phi}^{0}\in I,~\text{for all}~ u\in U.\]
    Assume $\|f^{u_{0}}\|_{\Phi}^{0}=0$, then $f^{u_{0}}=0$ and note that $\lambda\Big(\Phi\Big(\frac{|f|}{a}\Big)\Big)(u_{0})=0$, for all $a>0$. So for $\epsilon>0$, by continuity of $\lambda\Big(\Phi\Big(\frac{2|f|}{\epsilon}\Big)\Big)$, there exist open neighbourhood $U$ of $u_{0}$ such that $\lambda\Big(\Phi\Big(\frac{2|f|}{\epsilon}\Big)\Big)(u)<1$ which implies $\|f^{u}\|_{\Phi}^{0}<\epsilon$, for all $u\in U$. Hence proved.
\end{proof}
    Let $L^{\Phi}_{0}=\sqcup_{u\in G^{0}}L^{\Phi}(G^{u})$. Any $f \in C_{c}(G)$ can be identified as section $t_{f}:G^{0} \to L^{\Phi}_{0}$ where $t_{f}(u)= f_{|_{G^{u}}}$. Also for $g \in C_{c}(G^{u})$, by partition of unity argument, we can find $g’$ in $C_{c}(G)$, such that $g'_{|_{G^{u}}}=g$. These sections satisfy axioms (i), (ii) and (iii) of \cite[ Definition $10.1.2$]{MR0458185} and by  \cite[ Proposition $10.2.3$]{MR0458185}, $\left(\{L^{\Phi}(G^u)\}_{u\in G^{0}}, \Delta_{1}\right)$ becomes a continuous field of Orlicz space with continuous sections being sections locally close to $C_{c}(G)$, denoted by $\Delta_{1}$. Thus we can define a topology on the total space $L^{\Phi}_{0}=\sqcup_{u\in G^{0}}L^{\Phi}(G^{u})$ generated by  the following subsets: 
    \[U(V,\eta,\epsilon)=\{h\in L^{\Phi}_{0}:\|h-\eta(p(h))\|_{\Phi}^{0}<\epsilon, p(h)\in V\},\]
    where $\eta\in \Delta_{1}$ and $p:L^{\Phi}_{0}\to G^{0}$, is the projection of the total space on the base  $G^{0}$. We can also denote this continuous field of Orlicz space as $(L^{\Phi}_{0},\Delta_{1})$. 
    Denote the set of bounded sections of $L_{0}^{\Phi}$ as $I_{0}^{\Phi}(G, \lambda)$ and continuous sections vanishing at infinity as $E_{0}^{\Phi}$. Both are Banach spaces under the following norm:
\[\|\xi\|_{\Phi}^{0}=\sup _{u \in G^{0}}\left\|\xi^{u}\right\|_{\Phi}^{0}, \quad \xi^{u} \in L^{\Phi}\left(G^{u}\right). \] 

Now we have the following corollary. As the corollary follows similar lines to Theorem $3.1$, we skip it.
    \begin{corl}
    The family $\{M^\Psi(G^{u})\}_{u\in G^{0}}$ forms a continuous field of Morse-Transue space $(M^{\Psi}_{0}, \Delta )$ with $\Delta$ being the continuous sections, under $\|\cdot\|^{0}_{\Psi}$-norm,  locally close to $C_{c}(G)$ and $M^{\Psi}_{0}$ being the total space. $E^{\Psi}_{0}$ forms the continuous section  vanishing at infinity. 
\end{corl}
    
Note that $E_{0}^{\Phi}$ and $E_{0}^{\Psi}$ are $C_{b}\left(G^{0}\right)$ and $C_{0}\left(G^{0}\right)$-module: for $\xi \in E_{0}^{\Phi}, b \in C_{0}\left(G^{0}\right)$ or $C_{b}\left(G^{0}\right)$ we set $(\xi b)(u)=b(u) \xi^{u}$. 

The following proposition is the Orlicz space version of \cite[Proposition 7]{paterson2004fourier} and the proof is in similar lines.
\begin{prop}
     $C_{c}(G)$ is dense in $E_{0}^{\Phi}$ and $E_{0}^{\Psi}$.
\end{prop}
Now, we prove an important lemma which will be used for various results.
\begin{lem}
  Given $g \in L_{0}^{\Phi}\left(G^{u}\right)$ for any $u \in G^{0}$ with $\|g\|_{\Phi}^{0} \leq k$, there exist $\eta \in E_{0}^{\Phi}$ such that $\eta^{u}=g$ and $\|\eta\|_{\Phi}^{0} \leq k$.   
\end{lem}
\begin{proof}
    Assume $r=\|g\|_{\Phi}^{0} \neq 0$. If $r=0$, then the zero section is the required one. There exist a sequence $\left\{f_{n}\right\} \in C_{c}\left(G^{u}\right)$, tending to $g$, such that $\left\|f_{n}\right\|_{\Phi}^{0} \leq k, \forall n \in \mathbb{N}$. By passing to a subsequence if necessary we may suppose that $\left\|f_{n+1}-f_{n}\right\|_{\Phi}^{0}<\frac{k}{2^{n}}$. Then there exist $\left\{g_{n}\right\} \in C_{c}(G)$ such that $g_{0}^{u}=f_{1}, g_{n}^{u}=f_{n+1}-f_{n}$, and $\left\|g_{n}\right\|_{\Phi}^{0}<k 2^{-n}, n \geq 1$. The section $\eta_{1}=\sum_{n=0}^{\infty} g_{n}$ is in $E_{0}^{\Phi}$ and $\eta_{1}^{u}=\lim f_{n}=g$. By continuity of $\eta_{1}$, there exist an open neighbourhood $U$ around $u$ such that $\left\|\eta_{1}^{v}\right\|_{\Phi}^{0}>0, v \in U$. Take $b \in C_{c}\left(G^{0}\right)$ such that $b(u)=1,0 \leq b \leq 1$ and support of $b$ is in $U$. Then $\eta=r \eta_{1} \frac{b}{\left\|\eta_{1}\right\|_{\Phi}^{0}}$ is the required section. Here $\eta^{v}=r \eta_{1}^{v} \frac{b(v)}{\left\|\eta_{1}^{v}\right\|_{\Phi}^{0}},v\in G^{0}.$
    \end{proof}
    \begin{thm}
    Let $(\Phi, \Psi)$ be a complementary pair of $N$-function and $\Phi \in \Delta_{2}$. For $f \in C_{c}(G)$, the function $g: G^{0} \rightarrow \mathbb{R}$ defined as $g(u)=\left\|f^{u}\right\|_{\Phi}$ is a continuous function on $G^{0}$ with compact support. 
\end{thm}
\begin{proof}
     Let $u_{0} \in G^{0}, \epsilon>0$. Since $L^{\Phi}\left(G^{u_{0}}\right)=M^{\Psi}\left(G^{u_{0}}\right)^{*}$ and using  \cite[Chapter $3$, Theorem $13$]{MR1113700}, we can say that,
     \begin{align*}
        \left\|f^{u_{0}}\right\|_{\Phi}&=\sup\left\{\int_{G^{u_{0}}}|f g|d\lambda^{u_{0}}: \left\|g^{u_{0}}\right\|_{\Psi}^{0} \leq 1, g \in M^{\Psi}(G^{u_{0}}) \right\} \\ & =\inf \left\{\frac{1}{k}\left(1+\int_{G^{u_{0}}} \Phi(k|f|) d \lambda^{u_{0}}\right): k>0\right\}. 
     \end{align*}

By density of $C_{c}\left(G^{u_{0}}\right)$, there exist a $g^{\prime} \in C_{c}\left(G^{u_{0}}\right)$, such that $\left\|f^{u_{0}}\right\|_{\Phi}-\epsilon<\int\left|f g^{\prime}\right| d \lambda^{u_{0}}$. Hence, using partition of unity, we can find $g \in C_{c}(G)$ such that $g_{\left.\right|_{G} u_{0}}=g^{\prime}$ and $\|g\|_{\Psi}^{0} \leq 1$. By continuity of $\int|f g| d \lambda^{u}$, we can find a neighbourhood $U$ of $u_{0}$ such that $\left\|f^{u_{0}}\right\|_{\Phi}-\epsilon<\int|f g| d \lambda^{u}$ for all $u \in U$ which implies $\left\|f^{u_{0}}\right\|_{\Phi}-\epsilon<\left\|f^{u}\right\|_{\Phi}, \forall u \in U$.

By definition, $\exists$ $k>0$ such that $\frac{1}{k}\left(1+\int_{G^{u_{0}}} \Phi(k|f|) d \lambda^{u_{0}}\right)<\left\|f^{u_{0}}\right\|_{\Phi}+\epsilon$. So, by continuity, there exist an open neighbourhood $V$ of $u_{0}$, such that $$\frac{1}{k}\left(1+\int_{G^{u}} \Phi(k|f|) d \lambda^{u}\right)<\left\|f^{u_{0}}\right\|_{\Phi}+\epsilon, \forall u \in V.$$ Thus, $\forall u \in V \cap U,$ $$\left\|f^{u}\right\|_{\Phi} \in I=\left(\left\|f^{u_{0}}\right\|_{\Phi}-\epsilon,\left\|f^{u_{0}}\right\|_{\Phi}+\epsilon\right).$$ Hence proved.
\end{proof}
Here also we see $C_{c}(G)$ as sections satisfy axioms $(i),(ii)$ and $(iii)$ of \cite [Definition $10.1.2$]{MR0458185}, under $\| \cdot \|_{\Phi}$-norm and by \cite[ Proposition $10.2.3$]{MR0458185}, $\left(\left\{L^{\Phi}\left(G^{u}\right)\right\}_{u \in G^{0}}, \Delta_{1}^{\prime}\right)$ becomes a continuous field of Orlicz space with continuous sections being sections locally close to $C_{c}(G)$, denoted by $\Delta_{1}^{\prime}$. The continuous section vanishing at infinity is denoted by $E^{\Phi}$. It is a $C_{b}\left(G^{0}\right)$ and $C_{0}\left(G^{0}\right)$-module.

We can define a topology on the total space $L^{\Phi}=\sqcup_{u \in G^{0}} L^{\Phi}\left(G^{u}\right)$ generated by the following subsets:

\[
U(V, \eta, \epsilon)=\left\{h \in L^{\Phi}:\|h-\eta(p(h))\|_{\Phi}<\epsilon, p(h) \in V\right\},
\]
where $\eta \in \Delta_{1}^{\prime}$ and $p: L^{\Phi} \rightarrow G^{0}$, is the projection of the total space on the base $G^{0}$. We denote this continuous field of Orlicz space as $\left(L^{\Phi}, \Delta_{1}^{\prime}\right)$ also.

With the same idea as above, we can prove the following corollary.
\begin{corl}
     The family $\left\{M^{\Psi}\left(G^{u}\right)\right\}_{u \in G^{0}}$ becomes a continuous field of MorseTransue space $\left(M^{\Psi}, \Delta^{\prime}\right)$ with $\Delta^{\prime}$ being the continuous sections, under $\|\cdot\|_{\Psi}$-norm, locally close to $C_{c}(G)$ and $M^{\Psi}$ being the total space. $E^{\Psi}$ forms the continuous section vanishing at infinity. It is a $C_{b}\left(G^{0}\right)$ and $C_{0}\left(G^{0}\right)$-module.
\end{corl}
Using \cite[Lemma $2.3$]{bos2011continuous} , the continuous field of Orlicz spaces $\left(L_{0}^{\Phi}, \Delta_{1}\right)$ and $\left(L^{\Phi}, \Delta_{1}^{\prime}\right)$ are isomorphic. Here we can take the identity map as the morphism between two continuous fields of Banach spaces. Similarly, $\left(M^{\Psi}, \Delta^{\prime}\right)$ and $\left(M_{0}^{\Psi}, \Delta\right)$ are isomorphic. Note that Lemma $3.4$ can also be proved in the context of continuous field of Banach spaces $\left(L^{\Phi}, \Delta_{1}^{\prime}\right)$,$\left(M^{\Psi}, \Delta^{\prime}\right)$ and $\left(M_{0}^{\Psi}, \Delta\right)$. We can easily see that $C_{c}(G)$ is dense in $E^{\Phi}$ and $E^{\Psi}$.
\begin{defi}
    A strongly continuous representation of groupoids on a continuous field of Banach space $\left(B_{\pi}=\left\{B_{u}\right\}_{u \in G^{0}}, \Gamma\right)$ over $G^{0}$ is a triple $\left(B_{\pi}, \Gamma, \pi\right)$ such that
\begin{enumerate}[(i)]
    \item  $\pi(x) \in \mathcal{B}\left(B_{d(x)}, B_{r(x)}\right)$ is a invertible isometry, for each $x \in G$,

\item $\pi(u)$ is the identity map on $B_{u}$ for all $u \in G^{0}$,

\item $\pi(x) \pi(y)=\pi(x y)$ for all $(x, y) \in G^{2}$,

\item $\pi(x)^{-1}=\pi\left(x^{-1}\right)$ for all $x \in G$,

\item $x \rightarrow \pi(x) \eta^{d(x)}$ is continuous for every $\eta \in \Gamma$.
\end{enumerate}

\end{defi}
 \begin{thm}
     The left-regular representation $\left(L_{0}^{\Phi}, \Delta_{1}, L\right)$ is a strongly continuous representation. 
 \end{thm}
\begin{proof}
     Let $x \in G, r(x)=u, d(x)=v$ and $f \in L^{\Phi}\left(G^{v}\right)$. Then, $\int \Phi(|L(x) f|) d \lambda^{r(x)}=$ $\int L(x) \Phi(|f|) d \lambda^{r(x)}=\int \Phi(|f|) d \lambda^{d(x)}$. Thus it is an isometry. The conditions  $(i),(ii),$\\$(iii)$ and $(iv)$ of Definition $3.7$ can be easily verified.

We need to prove that $x \rightarrow L(x) \eta^{d(x)}$ is continuous for every $\eta \in \Delta_1$. Let $x_{0} \in G$,\\
then $L\left(x_{0}\right) \eta^{d\left(x_{0}\right)} \in L^{\Phi}\left(G^{r\left(x_{0}\right)}\right)$. For a given open set $K$ containing $L\left(x_{0}\right) \eta^{d\left(x_{0}\right)}$ in $L_{0}^{\Phi}$, by Lemma 3.4, we can find an open set of the form $U(V, \xi, \epsilon) \subset K$, where $V$ is a neighbourhood of $r\left(x_{0}\right), \epsilon>0$ and $\xi \in E_{0}^{\Phi}$ such that $L\left(x_{0}\right) \eta^{d\left(x_{0}\right)}=\xi^{r\left(x_{0}\right)}$. There exist $\eta_{0}, \xi_{0} \in C_{c}(G), V_{0} \subset V, r\left(x_{0}\right) \in V_{0}$ and an open set $V_{1}$ containing $d\left(x_{0}\right)$ such that:
\[
\left\|\xi^{w}-\xi_{0}^{w}\right\|_{\Phi}^{0}<\frac{\epsilon}{3}, \forall w \in V_{0},~\left\|\eta^{u}-\eta_{0}^{u}\right\|_{\Phi}^{0}<\frac{\epsilon}{3}, \forall u \in V_{1},~ L\left(x_{0}\right) \eta_{0}^{d(x_{0})}=\xi_{0}^{ r\left(x_{0}\right)} .
\]

Define $f \in C_{c}\left(G^{(2)}\right)$ as $f\left(x^{\prime}, y^{\prime}\right)=\Phi(\frac{3}{\epsilon}|\eta_{0}\left(x'^{-1}y^{\prime})-\xi_{0}(y^{\prime}\right)|)$. Using \cite[Lemma $3.12$]{bos2011continuous}, we can say that $F\left(x^{\prime}\right)=\int_{G^{r\left(x^{\prime}\right)}} \Phi\left(\frac{3}{\epsilon}|\eta_{0}\left(x'^{-1}y^{\prime}\right)-\xi_{0}\left(y^{\prime}\right)\right|) d \lambda^{r\left(x^{\prime}\right)}\left(y^{\prime}\right)$ is continuous on $G$. Note that $F\left(x_{0}\right)=0$ and by continuity of $F$, there exist an open neighbourhood $W$ of $x_{0}$ such that $F\left(x^{\prime}\right)<1, \forall x^{\prime} \in W$ which implies $$\left\|L\left(x^{\prime}\right) \eta_{0}^{d\left(x^{\prime}\right)}-\xi_{0}^{ r\left(x^{\prime}\right)}\right\|_{\Phi}^{0}\leq\frac{\epsilon}{3}$$ for every $x^{\prime} \in W$. So for all $x \in W^{\prime}=r^{-1}\left(V_{0}\right) \cap d^{-1}\left(V_{1}\right) \cap W$,
\begin{align*}
\left\|L(x) \eta^{d(x)}-\xi^{r(x)}\right\|_{\Phi}^{0} & \leq\left\|L(x) \eta^{d(x)}-L(x) \eta_{0}^{ d(x)}\right\|_{\Phi}^{0} \\
& +\left\|L(x) \eta_{0}^{d(x)}-\xi_{0}^{r(x)}\right\|_{\Phi}^{0}+\left\|\xi_{0}^{r(x)}-\xi^{r(x)}\right\|_{\Phi}^{0} \\
& =\left\|\eta^{d(x)}-\eta_{0}^{d(x)}\right\|_{\Phi}^{0}+\left\|L(x) \eta_{0}^{d(x)}-\xi_{0}^{r(x)}\right\|_{\Phi}^{0} \\
& +\left\|\xi_{0}^{r(x)}-\xi^{r(x)}\right\|_{\Phi}^{0}<\epsilon.
\end{align*}
Thus we have showed that, $L(x) \eta^{d(x)} \in U(V, \xi, \epsilon), \forall x \in W^{\prime}$. Hence we have proved that $x \rightarrow L(x) \eta^{d(x)}$ is continuous at $x_{0}$.
\end{proof}
Let $I(G, \lambda)$ denote the space of bounded sections of the Banach bundle $\left\{L^{1}(G)=\right.$ $\left.L^{1}\left(G^{u}, \lambda^{u}\right)\right\}_{u \in G^{0}}$ and $E^{1}$ denote the continuous sections over $G^{0}$ vanishing at infinity. Both form Banach spaces under the following norm:

\[
\|\eta\|_{1}=\sup _{u \in G^{0}}\left\|\eta^{u}\right\|_{1}, ~\eta^{u}\in L^{1}(G^{u}).
\]

It can be easily seen that $\left(L^{1}(G), \Gamma'\right)$ is a continuous field of Banach space where $\Gamma'$ is the space of continuous sections locally close to $C_{c}(G)$ under $\|\cdot\|_{1}$-norm. Note that $(C_{c}(G),\|\cdot\|_1)$ is a normed algebra under the following convolution:
\[f * g(x)  =\int_{G^{r(x)}} f(y) g\left(y^{-1} x\right) d \lambda^{r(x)}(y).\]
The density of $C_{c}(G)$ in $E^1$ by identifying $C_{c}(G)$ as sections, makes the space $E^1$ a Banach algebra.
\begin{lem}
     Let $G$ be a groupoid and $(\Phi,\Psi)$ be a complementary pair of $N$-functions with $\Phi \in\Delta_{2}$. Then $I_{0}^{\Phi}(G, \lambda) \subseteq I(G, \lambda)$, iff there exists $d>0$ such that $\|f\|_{1} \leq d\|f\|_{\Phi}^{0}, \forall f \in I_{0}^{\Phi}(G, \lambda)$.
\end{lem}
The proof is similar to that of \cite[lemma $3.3$]{kumar2019orlicz} using the Open Mapping Theorem. Similarly we can see that $E_{0}^{\Phi}\subseteq E^{1}$ iff  there exists $d>0$ such that $\|f\|_{1} \leq d\|f\|_{\Phi}^{0}, \forall f \in E_{0}^{\Phi}$.
\begin{corl}
  Let $(\Phi,\Psi)$ be a complementary pair of $N$-functions  with $\Phi \in \Delta_{2}$. If $I_{0}^{\Phi}(G, \lambda) \subseteq I(G, \lambda)$, then $E_{0}^{\Phi} \subseteq E^{1}$. 
\end{corl}
\begin{lem}
   Let $(\Phi,\Psi)$ be a complementary pair of $N$-functions  with $\Phi \in \Delta_{2}$. If $G$ is such that, $\sup _{u \in G^{0}} \lambda^{u}(G)<$ $\infty$, then $I_{0}^{\Phi}(G, \lambda) \subseteq I(G, \lambda)$ and $E^{\Phi}_{0}\subseteq E^{1}$. In particular, the conclusion holds if $G$ is a compact groupoid. 
\end{lem}
\begin{proof}

 Assume  that $\sup _{u \in G^{0}} \lambda^{u}(G)<\infty$. Since $\Phi$ is convex, there exists $c, r_{0}>0$ such that $\Phi(r) \geq c r, \forall r \geq r_{0} .$~$ \xi \in I_{0}^{\Phi}(G,\lambda)$ implies $ \sup _{u \in G^{0}} \int \Phi(a|\xi^{u}|) d \lambda^{u}<\infty$ for some $a>0$.

For $u \in G^{0}$, let $N^{u}=\left\{s \in G^{u}: a|\xi^{u}(s)|\leq r_{0}\right\}$.
\begin{align*}
\int|\xi^{u}| d \lambda^{u} & =\frac{1}{a}\left(\int_{N^{u}} a|\xi^{u}| d \lambda^{u}+\int_{G / N^{u}} a|\xi^{u}| d \lambda^{u}\right) \\
& \leq \frac{1}{a}\left(r_{0} \lambda^{u}(G)+\int_{G / N^{u}} a|\xi^{u}| d \lambda^{u}\right) \\
& \leq \frac{1}{a} \sup _{u \in G^{0}}\left(r_{0} \lambda^{u}(G)+\frac{1}{c} \int \Phi(a|\xi^{u}|) d \lambda^{u}\right) .
\end{align*}
Hence $\sup _{u \in G^{0}} \int|\xi^{u}| d \lambda^{u}<\infty $ implies $\xi \in I(G, \lambda)$.
\end{proof}

We will now show a sufficient condition for $E_{0}^{\Phi}$ to form a Banach algebra. It is the groupoid version of \cite[Theorem $3.5$]{kumar2019orlicz}. A necessary and sufficient condition is proved in \cite[Theorem $2$]{hudzik1987some}, when $G$ is a locally compact abelian group.

\begin{thm}
    Let $G$ be a locally compact groupoid and $(\Phi, \Psi)$ be a complementary pair of $N$-functions with $\Phi \in \Delta_{2}$. If $E_{0}^{\Phi} \subseteq E^{1}$, then the bilinear map $F: C_{c}(G) \times$ $C_{c}(G) \rightarrow C_{c}(G), F(f, g)=f * g$ is a bounded map under $\|\cdot\|^{0}_{\Phi}$-norm and can be extended to $E_{0}^{\Phi} \times E_{0}^{\Phi}$. In particular, the convolution is commutative if $G$ is a commutative group bundle.
\end{thm}
\begin{proof}
   Since $E_{0}^{\Phi} \subset E^{1}$, by Lemma $3.9$, there exists $d>0$ such that $\|f\|_{1} \leq d\|f\|^{0}_{\Phi}$, for all $f \in E_{0}^{\Phi}$. Let $f, g \in C_{c}(G)$ then $f * g(x)=\int f(y) g\left(y^{-1} x\right) d \lambda^{r(x)}(y) \in C_{c}(G)$.

Let $u \in G^{0}$. Using Fubini's theorem, Holder's Inequality, Lemma $3.4$ in the context of $M^{\Psi}_{0}$ and isometry of left translation $L(y): L^{\Phi}(G^{d(y)}) \rightarrow L^{\Phi}(G^{r(y)})$, we get
\begin{align*}
\left\|(f * g)^{u}\right\|_{\Phi} & =\sup \left\{\int\left|f * g\left\|h \mid d \lambda^{u}:\right\| h \|_{\Psi}^{0} \leq 1, h \in E_{0}^{\Psi}\right\}\right. \\
& \leq \sup \left\{\iint\left|f(y)\left\|g\left(y^{-1} x\right)\right\| h(x)\right| d \lambda^{r(x)}(y) d \lambda^{u}(x):\|h\|_{\Psi}^{0} \leq 1\right\} \\
& =\sup \left\{\iint\left|f(y)\left\|g\left(y^{-1} x\right)\right\| h(x)\right| d \lambda^{r(y)}(x) d \lambda^{u}(y):\|h\|_{\Psi}^{0} \leq 1\right\} \\
& \leq 2\|f\|_{1}\|g\|_{\Phi}^{0} .
\end{align*}
Thus,
\begin{align*}
\|(f * g)\|_{\Phi}^{0}=\sup _{u \in G^{0}}\left\|(f * g)^{u}\right\|_{\Phi}^{0} & \leq \sup _{u \in G^{0}}\left\|(f * g)^{u}\right\|_{\Phi} 
 \leq 2\|f\|_{1}\|g\|_{\Phi}^{0}. 
\end{align*}
So, $\|F(f, g)\|_{\Phi}^{0} \leq 2 d\|f\|_{\Phi}^{0}\|g\|_{\Phi}^{0}$. Now by density of $C_{c}(G)$, we extend $F$ to $E_{0}^{\Phi} \times E_{0}^{\Phi}$ and $\|F(\xi,\eta)\|_{\Phi}^{0} \leq 2 d\|\xi\|_{\Phi}^{0}\|\eta\|_{\Phi}^{0}$. Define an equivalent norm $\|\cdot\|_{\Phi}^{\prime}=2 d\|\cdot\|_{\Phi}^{0}$, so that $\|F(f, g)\|_{\Phi}^{\prime} \leq\|f\|_{\Phi}^{\prime}\|g\|_{\Phi}^{\prime}$.

If $G$ is a commutative group bundle, $G^{u}=G_{u}=G(u)$ and $\lambda^{u}=\lambda_{u}$ for all $u \in G^{0}$. Hence for $f, g \in C_{c}(G)$,
\begin{align*}
f * g(x) & =\int_{G^{r(x)}} f(y) g\left(y^{-1} x\right) d \lambda^{r(x)}(y) \\
& =\int_{G_{r(x)}} f\left(y^{-1}\right) g(y x) d \lambda_{r(x)}(y) \\
& =\int_{G_{d(x)}} \check{f}\left(y x^{-1}\right) g(y) d \lambda_{d(x)}(y) \\
& =\int_{G_{d(x)}} f\left(x y^{-1}\right) g(y) d \lambda_{d(x)}(y) \\
& =\int_{G^{r(x)}} f\left(y^{-1} x\right) g(y) d \lambda^{r(x)}(y)=g*f(x). 
\end{align*}
\end{proof}
For $\xi, \eta \in E_{0}^{\Phi}$, we denote $F(\xi, \eta)$ as $\xi * \eta$. We showed the sufficient condition for $E_{0}^{\Phi}$ to become a Banach algebra, by taking the equivalent norm mentioned above. We denote it again as $\|\cdot\|_{\Phi}^{0}$ for convenience. Similarly, we can see that algebra structure exists for $E^{\Phi}$ also under the same sufficient condition.

\begin{lem}
     Let $G$ be a groupoid and  $(\Phi, \Psi)$ be a complementary pair of $N$-function with $\Phi \in \Delta_{2}$. For $f \in C_{c}(G)$, the map $L_{f}: C_{c}(G) \rightarrow E_{0}^{\Phi}$, defined as $L_{f}(g)=f * g$, extends to a bounded linear map on $E_{0}^{\Phi}$ with   $\left\|L_{f}\right\| \leq\|f\|_{1}$. Moreover, the norm decreasing homomorphism $f\to L_f$ from $(C_{c}(G),\|\cdot\|_{1})$ into $B(E^{\Phi}_{0})$ can be extended to $E^{1}$.
\end{lem}
\begin{proof}
    Let $f , g \in C_{c}(G)$ with $\|f\|_{1}=1$ and $\|g\|_{\Phi}^{0}=1$. We will show that $\int \Phi(|f * g|)\ d\lambda^{u}<1, \forall u \in G^{0}$. Let $u \in G^{0}$ and $\left\|f^{u}\right\|_{1} \neq 0$, using Jensen's inequality, Fubini's theorem and isometry of left translation,
\begin{align*}
\int \Phi(|f * g|)(x) d \lambda^{u}(x) & =\int \Phi\left(\left|\int f(y) g\left(y^{-1} x\right) d \lambda^{r(x)}(y)\right|\right) d \lambda^{u}(x) \\
& \leq \int \Phi\left(\int|f(y)|\left|g\left(y^{-1} x\right)\right| d \lambda^{u}(y)\right) d \lambda^{u}(x) \\
& \leq \int \Phi\left(\int \frac{|f(y)|\left|g\left(y^{-1} x\right)\right|}{\left\|f^{u}\right\|_{1}} d \lambda^{u}(y)\right) d \lambda^{u}(x) \\
& \leq \iint \Phi\left(\left|g\left(y^{-1} x\right)\right|\right) d \nu^{u}(y) d \lambda^{u}(x) \\
& =\iint \Phi\left(\left|g\left(y^{-1} x\right)\right|\right) d \lambda^{u}(x) d \nu^{u}(y) 
\leq 1,
\end{align*}
where $\nu^{u}(E)=\frac{1}{\|f^{u}\|_{1}}\int |f|d\lambda^{u},$ $E \subset G^{u}.$ Thus $L_{f}$ can be extended to $E_{0}^{\Phi}$ and $\left\|L_{f}\right\| \leq\|f\|_{1}$. Also by density of $C_{c}(G)$ in $E^{1}$ we can extend the norm decreasing homomorphism $f \to L_f$ and thus define $L_{\xi}$ for every $\xi \in E^{1}$. 
\end{proof}
The above lemma shows that $E_{0}^{\Phi}$ is an $E^{1}$-module and denote $L_{\xi}(g)$ as $\xi * g$ for $\xi \in E^{1}$ and $g \in E_{0}^{\Phi}$. This is analogous to \cite [Corollary $3.4 (i)$]{aghababa2020fig}.

In the following lemma, we show that the right convolution of $C_{c}(G)$ over $E^{\Phi}_{0}$ is a bounded operator.
\begin{lem}
    Let $G$ be a groupoid and $(\Phi,\Psi)$ be a complementary pair of $N$- functions both satisfying $\Delta_{2}$-condition. For $F \in C_{c}(G)$, the map $R_{F}:C_{c}(G)\to C_{c}(G)$, defined as $R_{F}(g)= g*F$, extends to a bounded linear map on $E^{\Phi}_{0}$ with $\|R_{F}\|\leq 2 K_{F}^2$, where  $K_{F}= \max\{\sup_{u\in G^{0}}\|\Phi^{-1}\circ |F^{u}|\|_{\Tilde{\Psi}}^{0}, \sup_{u\in G^{0}}\|\Psi^{-1}\circ |\check{F}^{u}|\|_{\Tilde{\Psi}}^{0}\}$, for an $N$- function $\Tilde{\Psi}$ . 
\end{lem}
\begin{proof}
    Let $F,g, h\in C_{c}(G)$ and $ u \in G^{0}$.
 \begin{align*}
& \left| \int g*F(x)h(x)d\lambda^{u}(x)\right| \leq  \int \int |g(y)F(y^{-1}x)||h(x)|d\lambda^{r(x)}(y)d\lambda^{u}(x)\\
\leq & \int\int \left[|g(y)|\Phi^{-1}(|F(y^{-1}x)|)\right]\left[|h(x)|\Psi^{-1}(|F(y^{-1}x)|)\right]d\lambda^{u}(y)d\lambda^{u}(x) \leq 2 AB
\end{align*}
where
\begin{align*}
A &= \|k(x,y)\|_{\Phi}^{0},~ k(x,y)= |g(y)|\Phi^{-1}(|F(y^{-1}x)|)\\
B &= \|m(x,y)\|_{\Psi}^{0},~ m(x,y)= |h(x)|\Psi^{-1}(|F(y^{-1}x)|),~~ (x,y)\in G^{u}\times G^{u}. 
\end{align*}
According to \cite[Proposition $2.13$ and Proposition $2.4(ii)$]{MR1113700}, there exist an $N$- function $\Tilde{\Psi}= \max\{\Psi_{1},\Psi_{2}\}$ where $\Psi_{1}(a)= \sup\{\frac{\Phi(ab)}{\Phi(b)}: b\geq 0\}$ and $\Psi_{2}(a)=\sup\{\frac{\Psi(ab)}{\Psi(b)}: b\geq 0\}~ a,b\geq 0$, so that $\Phi(ba)\leq \Phi(b)\Tilde{\Psi}(a)$ and $\Psi(ba)\leq \Psi(b)\Tilde{\Psi}(a)$. Also note that the functions $\Phi^{-1}\circ |F|$ and $\Psi^{-1}\circ |\check{F}|$ are bounded and vanishes outside the supp$(F)$ and supp$(\check{F})$ respectively. So,
\begin{eqnarray*}
 & & \int\int \Psi\Bigg(\frac{|h(x)|\Psi^{-1}(|F(y^{-1}x)|)}{\|h\|^{0}_{\Psi} \|\Psi^{-1}\circ \check{F}\|_{\Tilde{\Psi}}^{0}}\Bigg)d\lambda^{u}(x)d\lambda^{u}(y)\\
 & \leq & \int\int \Psi\left(\frac{|h(x)|}{\|h\|^{0}_{\Psi}}\right)\Tilde{\Psi}\left(\frac{\Psi^{-1}(|F(y^{-1}x)|)}{\|\Psi^{-1}\circ \check{F}\|_{\Tilde{\Psi}}^{0}}\right)d\lambda^{u}(x)d\lambda^{u}(y)\\
 & = & \int \Psi\left(\frac{|h(x)|}{\|h\|^{0}_{\Psi}}\right) \int \Tilde{\Psi}\left( \frac{\Psi^{-1}(|\check{F}(y)|)}{\|\Psi^{-1}\circ \check{F}\|_{\Tilde{\Psi}}^{0}}\right)d\lambda^{d(x)}(y) d\lambda^{u}(x) \leq 1.
\end{eqnarray*}
Thus, 
\begin{equation*}
    \|m(x,y)\|^{0}_{\Psi} \leq \|h\|^{0}_{\Psi} \|\Psi^{-1}\circ |\check{F}|\|_{\Tilde{\Psi}}^{0}.
\end{equation*}
Similarly, 
\begin{equation*}
    \|k(x,y)\|^{0}_{\Phi} \leq \|g\|^{0}_{\Phi} \|\Phi^{-1}\circ |{F}|\|_{\Tilde{\Psi}}^{0}.
\end{equation*}
Hence, we can see that 
\begin{equation*}
    \|g*F\|_{\Phi}^{0}\leq 2 \|g\|^{0}_{\Phi} K_{F}^2.
\end{equation*}

\end{proof}

Now, we will show the existence of the left approximate identity on $E_{0}^{\Phi}$ when it is a Banach algebra. The corresponding results are proved in \cite[Proposition $2$]{rao2001convolutions} for the group case and \cite[Theorem $3.8$]{kumar2019orlicz}  for the hypergroup case.
\begin{thm}
     Let $G$ be a groupoid, and let $(\Phi, \Psi)$ be a complementary pair of $N$-functions with $\Phi \in \Delta_{2}$, then there exist a sequence $\{e_{n}\}$ of continuous functions with compact support  bounded in $\|\cdot\|_{1}$-norm such that  $\|e_{n}*\xi-\xi\|_{\Phi}^{0}\to 0$ for every $\xi \in E_{0}^{\Phi}$.
\end{thm}
\begin{proof}
     Let $f \in C_{c}(G)$ with $\operatorname{supp}(f)=K$ and $\epsilon>0$. By \cite[Proposition $11$]{paterson2004fourier}, there exists a bounded approximate identity $\left\{e_{n}\right\} \geq 0$ for $C_{c}(G)$ under $\|\cdot\|_{1}$-norm with $\left\|e_{n}\right\|_{1} \leq 2$, for all $n \in \mathbb{N}$. There exist a sequence $\left\{U_{n}\right\}$ of open neighbourhoods of $G^{0}$ in $G$ such that each $U_{n}$ is $d$-relatively compact, $U_{n} \subset U_{1}$ and is a fundamental sequence for $G^{0}$ in the sense that every neighbourhood $V$ of $G^{0}$ in $G$ contains $U_{n}$ eventually. There is an increasing sequence $\left\{K_{n}\right\}$ of compact subsets of $G^{0}$ such that $\cup K_{n}=G^{0}$. Also $\left\{e_{n}\right\}$ is constructed such that $\operatorname{supp}\left(e_{n}\right) \subset U_{n}$ and $\int e_{n} d \lambda^{u}=1$ for all $u \in K_{n}$.\\
Let $L=\overline{U_{1} K}$, where $U_{1} K=\left(U_{1} \cap d^{-1}(r(K))\right) K$ is relatively compact such that $\operatorname{supp}(f)$ and $\operatorname{supp}\left(e_{n} * f\right)$ is contained in $L$. Since $L$ is compact, $\sup _{u \in G^{0}} \lambda^{u}(L)<\infty$ and let that value be $M$. Choose $\delta>0$ such that $\Phi(x)<\frac{1}{2 M}$ for all $x$ with $|x|<\delta$. This is possible by continuity of $\Phi$. Note that $r(L) \subset K_{n}$ for all $n \geq n_{0}$ for some $n_{0} \in \mathbb{N}$. Also by continuity of $f$ and product in $G$, there exist $n_{1} \in \mathbb{N}, n_{1} \geq n_{0}$ such\\
that for $n \geq n_{1}, \frac{\left|f\left(y^{-1} x\right)-f(x)\right|}{\epsilon}<\delta$ for all $(x, y) \in\left(L \times U_{n}\right) \cap G^{2}$. Now choose $e_{n}$, where $n \geq n_{1}$. Let $u \in r(L)$,

\begin{align*}
\int \Phi\left(\frac{\left|e_{n} * f-f\right|}{\epsilon}\right)(x) d \lambda^{u}(x) & =\int_{L} \Phi\left(\frac{1}{\epsilon}\left|\int e_{n}(y) f\left(y^{-1} x\right) d \lambda^{r(x)}(y)-f(x)\right|\right) d \lambda^{u}(x) \\
& \leq \int_{L} \Phi\left(\int e_{n}(y)\frac{\left|f\left(y^{-1} x\right)-f(x)\right|}{\epsilon} d\lambda^{r(x)}(y)\right) d \lambda^{u}(x) \\
& =\int_{L} \Phi\left(\int \frac{\left|f\left(y^{-1} x\right)-f(x)\right|}{\epsilon} d\nu_{n}^{u}(y)\right) d \lambda^{u}(x)\\
& \quad \text { where } \nu_{n}^{u}(E)=\int_{E} e_{n}(y) d \lambda^{u}(y),~E\subset G^{u}\\
& \leq \int_{L} \int_{U_{n}} \Phi\left(\frac{\left|f\left(y^{-1} x\right)-f(x)\right|}{\epsilon}\right) d \nu_{n}^{u}(y) d \lambda^{u}(x) \\
& =\int_{U_{n}} \int_{L} \Phi\left(\frac{\left|f\left(y^{-1} x\right)-f(x)\right|}{\epsilon}\right) d \lambda^{u}(x) d \nu_{n}^{u}(y) \\
& \leq \frac{1}{2 M} M \int e_{n}(y) d \lambda^{u}\leq 1.
\end{align*}

The above is true for all $u \in r(L)$ and all $n \geq n_{1}$.
Also, for $u \notin r(L)$,\\$ \int \Phi\left(\frac{\left|e_{n} *-f\right|}{\epsilon}\right) d \lambda^{u}=0$, since $\operatorname{supp}\left(e_{n} * f-f\right) \subset L$.\\
So, $\left\|e_{n} * f-f\right\|_{\Phi}^{0} \leq \epsilon$ for all $n \geq n_{1}$. The result follows by Lemma $3.13$ and the density of $C_{c}(G)$.
\end{proof}
Using \cite [Theorem $4.6$]{lazar2018selection}, for a closed $C_{b}\left(G^{0}\right)$-submodule $I$ of $E_{0}^{\Phi}$, there exist a  subbundle of $(L_{0}^{\Phi},p,G^{0})$, denoted by $(I_{0}^{\Phi},p_{I_{0}^{\Phi}},G^{0})$, such that the set of all sections of $(I_{0}^{\Phi},p_{I_{0}^{\Phi}},G^{0})$ vanishing at infinity on $G^{0}$ coincides with $I$. Here the fiber $I_{0}^{\Phi}(u)=\left\{\xi^{u}: \xi \in I\right\}$ is a closed subspace of $L^{\Phi}\left(G^{u}\right)$. We say that the subbundle $(I_{0}^{\Phi},p_{I_{0}^{\Phi}},G^{0})$ is invariant under the groupoid representation $\pi$, if $\pi(x)(I_{0}^{\Phi}(d(x))) \subset I_{0}^{\Phi}(r(x))$ for all $x \in G$.

The next theorem provides the sufficient and necessary condition for a closed $C_{b}\left(G^{0}\right)$ submodule of $E_{0}^{\Phi}$ to become an ideal when $E_{0}^{\Phi}$ is a Banach algebra. This is the groupoid version of the result in \cite[Theorem $3.9$]{kumar2019orlicz}.
\begin{thm}
    Let $G$ be a groupoid, $(\Phi, \Psi)$ be a complementary pair of $N$-functions with $\Phi \in \Delta_{2}$ and $E_{0}^{\Phi}$ be an algebra. A closed $C_{b}\left(G^{0}\right)$-submodule $I$ of $E_{0}^{\Phi}$ is a left ideal if and only if the  subbundle $(I_{0}^{\Phi},p_{I_{0}^{\Phi}},G^{0})$  is invariant under left-regular representation of $G$.
\end{thm}
\begin{proof}
    Let  $u \in G^{0}, f, g \in C_{c}(G).$ $ L_{z}$ denote the left translation from $L^{\Phi}(G^{d(z)})$ to $L^{\Phi}(G^{r(z)})$. For $x\in G^{r(z)},$
\begin{align*}
L_{z}(f * g)^{d(z)}(x)  =(f * g)^{d(z)}\left(z^{-1} x\right)&=\int_{G^{d(z)}} f(y) g\left(y^{-1} z^{-1} x\right) d \lambda^{d(z)}(y) \\
& =\int_{G^{d(z)}} f(y) \check{g}\left(x^{-1} z y\right) d \lambda^{d(z)}(y) \\
& =\int_{G^{r(z)}} L_{z} f(y) g\left(y^{-1} x\right) d \lambda^{r(z)}(y) \\
& =\left(L_{z} f * g\right)(x),
\end{align*}
where $(L_{z} f * g)$ denotes the function $(\tilde{f}*g)^{r(z)}$ for $\tilde{f} \in C_{c}(G)$ such that $\tilde{f}_{|_{G^{r(z)}}}=L_zf^{d(z)}$.
So by density of $C_{c}(G),~ L_{z}(f * \eta)^{d(z)}=(\tilde{f} * \eta)^{r(z)}\in L^{\Phi}(G^{r(z)})$ for $\eta \in E_{0}^{\Phi}, \tilde{f} \in C_{c}(G)$ and can be denoted as $(L_{z} f * \eta)$.

Let $I$ be a closed left ideal of $E_{0}^{\Phi}$. Suppose $\left\{e_{n}\right\}$ is the left approximate identity of $E_{0}^{\Phi}$. Let $x \in G, d(x)=u, r(x)=v$ and $\eta^{u}\in I_{0}^{\phi}(u)$ for $\eta \in I$.
\begin{align*}
\left\|\left(L_{x} e_{n} * \eta-L_{x} \eta^{u}\right)\right\|_{\Phi}^{0} &= \left\|L_{x}\left(e_{n} * \eta-\eta\right)^{u}\right\|_{\Phi}^{0} \\
& =\left\|\left(e_{n} * \eta-\eta\right)^{u}\right\|_{\Phi}^{0} \leq\left\|\left(e_{n} * \eta-\eta\right)\right\|_{\Phi}^{0} \to 0.
\end{align*}

As mentioned before $\left(L_{x} e_{n} * \eta\right)$ denotes the function $\left(e_{n}^{\prime} * \eta\right)^{v}$ for some $e_{n}^{\prime} \in C_{c}(G)$ with ${e_n}'_{|_{G^{v}}}=L_{x}e_{n}^{u}$ and since $I$ is a left ideal, $e_{n}^{\prime} * \eta \in I$ for every $n \in \mathbb{N}$. Hence $\left(e_{n}^{\prime} * \eta\right)^{v} \in I_{0}^{\Phi}(v)$. So $L_{x} \eta^{u} \in I_{0}^{\Phi}(v)$. Since $x$ was arbitrary we can say that $(I_{0}^{\Phi},p_{I_{0}^{\Phi}},G^{0})$  is invariant under left regular representation.

Let $(I_{0}^{\Phi},p_{I_{0}^{\Phi}},G^{0})$ is invariant under left regular representation. If $I$ is not an ideal, by density of $C_{c}(G)$ in $E^{\Phi}_{0}$, there exist $f' \in C_{c}(G)$ and $g \in I$ such that $f' * g \notin I$. Hence by \cite [Theorem $4.6$]{lazar2018selection}, there exist $u_{0} \in G^{0}$, such that $(f' * g)^{u_{0}} \notin I_{0}^{\Phi}\left(u_{0}\right)$. So there exist $h \in L^{\Psi}\left(G^{u_{0}}\right)$ such that $\int_{G^{u_{0}}}\left(f^{\prime} * g\right) h d \lambda^{u_{0}} \neq 0$ and $\int_{G^{u_{0}}} k h d \lambda^{u_{0}}=0$, for all $k \in I_{0}^{\Phi}\left(u_{0}\right)$.\\
$f^{\prime} * g_{n} \rightarrow f^{\prime} * g$ in $E_{0}^{\Phi}$ as $g_{n} \rightarrow g, g_{n} \in C_{c}(G)$. So in particular, by Holder's inequality, $\int_{G^{u_{0}}}\left(f^{\prime} * g_{n}\right) h d \lambda^{u_{0}} \rightarrow \int_{G^{u_{0}}}\left(f^{\prime} * g\right) h d \lambda^{u_{0}}$.
\begin{align*}
\int_{G^{u_{0}}}\left(f^{\prime} * g_{n}\right) h d \lambda^{u_{0}} & =\int_{G^{u_{0}}}\left(\int_{G^{u_{0}}} f^{\prime}(y) g_{n}\left(y^{-1} x\right) d \lambda^{u_{0}}(y)\right) h(x) d \lambda^{u_{0}}(x) \\
& =\iint_{G^{u_{0}}} f^{\prime}(y) g_{n}(y^{-1} x) h(x) d \lambda^{u_{0}}(y) d \lambda^{u_{0}}(x) \\
& =\iint_{G^{u_{0}}} g_{n}(y^{-1} x) h(x) f^{\prime}(y) d \lambda^{u_{0}}(x)  d \lambda^{u_{0}}(y).
\end{align*}

Note that since $L_{y}:L^{\Phi}(G^{d(y)})\to L^{\Phi}\left(G^{u_{0}}\right)$ is isometry, the above integral converges to $\iint_{G^{u_{0}}} L_{y} g(x) h(x) f^{\prime}(y) d \lambda^{u_{0}}(x) d \lambda^{u_{0}}(y)$. Since $L_{y} g^{d(y)}\in I_{0}^{\Phi}\left(u_{0}\right)$,\\$\iint_{G^{u_{0}}} L_{y} g(x) h(x) f^{\prime}(y) d \lambda^{u_{0}}(x) d \lambda^{u_{0}}(y)=0$, for all $y \in G^{u_{0}}$.\\ So $\int_{G^{u_{0}}}\left(f^{\prime} * g\right) h d \lambda^{u_{0}}=0$, which is a contradiction. Hence proved.
\end{proof}
\section{The space of Convolutors of $E_{0}^{\Phi}$}
In locally compact groups and hypergroups, the space of convolutors of Morse-Transue space $M^{\Phi}(G)$, is identified with the dual of a space, denoted by $\check{A}_{\Phi}(G)$. This well-known result can be found in \cite{aghababa2013space} and \cite{kumar2019orlicz}. In this section, we provide a groupoid version of this result. Here we assume $(\Phi, \Psi)$ as complementary pair of $N$-functions both satisfying  $\Delta_{2}$-condition.

A bounded linear operator $T$ on $E_{0}^{\Phi}$ is called a convolutor if 
$T(f * g)=$ $Tf * g$, for all $f, g \in C_{c}(G)$. This space is denoted as $C_{V_{\Phi}}(G)$. It can be easily verified that $C_{V_{\Phi}}(G)$ is a closed subspace of $B(E_{0}^{\Phi})$.
\begin{lem}
 Let $G$ be a groupoid and $(\Phi, \Psi)$ be a complementary pair of $N$-functions. If $T \in C_{V_{\Phi}}(G)$, then there exist $\{e_{n}\} \in C_{c}(G)$ with $\|e_{n}\|_{1} \leq 2$ such that if we set $T_{n}(f)=T\left(e_{n} * f\right), f \in E^{\Phi}_{0}$,    
\begin{enumerate}[(i)]
\item $\left\|T_{n}\right\| \leq 2\|T\|$, 
\item  $\lim _{n \rightarrow\infty}\left\|T_{n} f-T f\right\|_{\Phi}^{0}=0, \text { for } f \in E^{\Phi}_{0}.$
\end{enumerate}
\end{lem}
\begin{proof}
    We proved that there exist $\left\{e_{n}\right\} \in C_{c}(G)$ with $\|e_{n}\|_{1} \leq 2$ and $\left\|e_{n} * f-f\right\|_{\Phi}^{0} \rightarrow 0$ for all $f \in E_{0}^{\Phi}$.
Take $f \in C_{c}(G)$,
\[
\left\|T\left(e_{n} * f\right)\right\|_{\Phi}^{0} \leq\|T\|\left\|e_{n}\right\|_{1}\|f\|_{\Phi}^{0} \leq 2\|T\|\|f\|_{\Phi}^{0} .
\]
We can see that $T_{n}$ is contained in $C_{V_{\Phi}}(G)$ and $\left\|T_{n}(f)\right\|_{\Phi}^{0} \leq 2\|T\|\|f\|_{\Phi}^{0}$ for all $f \in E_{0}^{\Phi}$. Thus, $\left\|T_{n}\right\| \leq 2\|T\|$ for every $n \in \mathbb{N}$. For $f \in C_{c}(G)$,
\begin{align*}
\lim _{n \rightarrow \infty}\left\|T_{n}(f)-T(f)\right\|_{\Phi}^{0}=&\lim _{n \rightarrow \infty}\left\|T\left(e_{n} * f-f\right)\right\|_{\Phi}^{0} \\
\leq&\|T\| \lim _{n \rightarrow \infty}\left\|e_{n} * f-f\right\|_{\Phi}^{0}=0.
\end{align*}
\end{proof}
 Let $\mathcal{K}(G)$ be the collection of compact sets with nonempty interior intersecting $G^{0}$. For $P \in \mathcal{K}(G)$, define 
\begin{equation*}
E^{\Phi}_{0}(P)=  \overline{\{f\in C_{c}(G): \operatorname{supp}(f)\subset P\}}^{\|\cdot\|^{0}_{\Phi}}.
\end{equation*}
\begin{equation*} 
E^{\Psi}_{0}(P)= \overline{\{f\in C_{c}(G): \operatorname{supp}(f)\subset P\}}^{\|\cdot\|^{0}_{\Psi}}.
\end{equation*}
By the density of $C_{c}(G)$ and Holder's inequality, we can define a function in $C_{c}(G)$ denoted as $\xi*\check{\eta}$, for $\xi \in E^{\Psi}_{0}(P)$ and $\eta \in E^{\Phi}_{0}(P)$.
\begin{align*}
 \check{A}_{\Phi, P}(G)=\Big\{h \in C_{c}(G): h=\sum_{n=1}^{\infty} g_{n} * \check{f}_{n}, f_{n} \in E^{\Phi}_{0}(P)&, ~g_{n} \in E^{\Psi}_{0}(P), \\
&\sum_{n=1}^{\infty}\left\|f_{n}\right\|_{\Phi}^{0}\left\|g_{n}\right\|_{\Psi}^{0}<\infty\Big\} .
\end{align*}
The norm is defined as follows,
\[
\|h\|_{\check{A}_{\Phi, P}(G)}=\inf \left\{\sum_{n=1}^{\infty}\left\|f_{n}\right\|_{\Phi}^{0}\left\|g_{n}\right\|_{\Psi}^{0}: h=\sum_{n=1}^{\infty} g_{n} * \check{f}_{n}\right\}.
\]
Set $\check{A}_{\Phi}(G)=\cup_{P \in \mathcal{K}(G)} \check{A}_{\Phi, P}(G)$. Then $\check{A}_{\Phi}(G)$ is a subspace of $C_{c}(G)$ and for $h \in$ $\check{A}_{\Phi}(G)$,
\[
\|h\|_{\check{A}_{\Phi}(G)}=\inf \left\{\|h\|_{\check{A}_{\Phi, P}(G)}: h \in \check{A}_{\Phi, P}(G), P \in \mathcal{K}(G)\right\}.
\]

For $f, g \in C_{c}(G)$, since $|f * \check{g}(x)| \leq 2\|f\|_{\Phi}^{0}\|g\|_{\Psi}^{0}$, for all $x \in G,\|h\|_{\infty} \leq 2\|h\|_{\check{A}_{\Phi}(G)}$. $\bar{A}_{\Phi}(G)$ denotes the norm completion of $\check{A}_{\Phi}(G)$. When $G$ is a locally compact group, this space coincides with the norm completion of $\check{A}_{\Phi}(G)$, defined in \cite[Section $4$, Page $27$]{aghababa2013space}.
\begin{lem}
  $\check{A}_{\Phi}(G)$ is a left and right $C_{0}\left(G^{0}\right)$-module with following action:

For $b \in C_{0}(G^{0}), h \in \check{A}_{\Phi}(G),(bh)(x)=b(d(x)) h(x)$ and $(h b)(x)=h(x) b(r(x))$.   
\end{lem}
\begin{proof}
     Let $f, g \in C_{c}(G)$, having support on $P \in \mathcal{K}(G)$, then $h=f * \check{g} \in \check{A}_{\Phi}(G)$.
\begin{align*}
(h b)(x) & =b(r(x)) h(x)=b(r(x)) f * \check{g}(x) \\
& =b(r(x)) \int f(y) \check{g}\left(y^{-1} x\right) d \lambda^{r(x)}(y) \\
& =\int f(y) \check{g}\left(y^{-1} x\right) b(r(y)) d \lambda^{r(x)}(y) \\
& =f b * \check{g}(x) \in \check{A}_{\Phi}(G) .
\end{align*}

Also $(f b) \in E^{\Psi}_{0} (P)$. Hence $h b \in \check{A}_{\Phi}(G)$ for any $h \in \check{A}_{\Phi}(G)$. Similarly, we can see that $b h(x)=f *\check{(g b)}(x)$. Here, $\|h b\|_{\check{A}_{\Phi}(G)} \leq\|b\|_{\infty}\|h\|_{\check{A}_{\Phi}(G)}$ and $\|b h\|_{\check{A}_{\Phi}(G)} \leq$ $\|b\|_{\infty}\|h\|_{\check{A}_{\Phi}(G)}$.
\end{proof}
 Let $B_{D}^{\Phi}\left(\check{A}_{\Phi}(G), C_{0}(G^{0})\right)$ be the space of all bounded right-module linear maps from $\check{A}_{\Phi}(G)$ to $C_{0}(G^{0})$. For $\alpha \in B_{D}^{\Phi}\left(\check{A}_{\Phi}(G), C_{0}(G^{0})\right)$, $f \in C_{c}(G)$, \\define $f \alpha: C_{c}(G) \rightarrow C_{0}(G^{0})$ as, $f \alpha(g)=\alpha(g * \check{f})$. It is a linear map. Also,
\[
|f \alpha(g)(u)|=|\alpha(g * \check{f})(u)| \leq\|\alpha\|\|g * \check{f}\|_{\check{A}_{\Phi}(G)} \leq\|\alpha\|\|g\|_{\Psi}^{0}\|f\|_{\Phi}^{0} \leq\|\alpha\|\|g\|_{\Psi}\|f\|_{\Phi}^{0}.
\]

Hence, $\|f \alpha(g)\|_{\infty} \leq\|\alpha\|\|f\|_{\Phi}^{0}\|g\|_{\Psi}$. We can extend $f \alpha$ to $E^{\Psi}$.

$f \alpha(g b)=\alpha(g b * \check{f})=\alpha((g * \check{f}) b)=\alpha((g * \check{f})) b$.
Thus, $\alpha \rightarrow f \alpha$ is a bounded linear map from $B_{D}^{\Phi}\left(\check{A}_{\Phi}(G), C_{0}(G^{0})\right)$ to $B_{D}^{\Phi}\left(E^{\Psi}, C_{0}(G^{0})\right)$, where $B_{D}^{\Phi}\left(E^{\Psi}, C_{0}(G^{0})\right)$ is the space of bounded linear right-module maps from $E^{\Psi}$ to $C_{0}(G^{0})$.\\
For $\xi \in E^{\Phi}_{0}$, $\eta \in E^{\Psi}$, define:
\[\langle\xi,\eta\rangle(u)=\int_{G^{u}}\xi^{u}\eta^{u}d\lambda^{u}.\]

 Note that the function $\langle \xi, \eta\rangle$ is in $C_{0}\left(G^{0}\right)$ such that, $$\|\langle \xi, \eta \rangle\|_{\infty} \leq \|\xi\|_{\Phi}^{0}\|\eta\|_{\Psi}\leq 2\|\xi\|_{\Phi}^{0}\|\eta\|_{\Psi}^{0}.$$ This is possible due to denseness of $C_{c}(G)$ and Holder's inequality.\\
Let, \[R\left(E^{\Psi}, C_{0}(G^{0})\right)=\left\{S \in B_{D}\left(E^{\Psi}, C_{0}(G^{0})\right): S(\xi)(u)=\langle g, \xi\rangle(u) \, \text{for some} \, g \in E_{0}^{\Phi}(G)\right\}.\] Since $M^{\Psi}\left(G^{u}\right)^{*}=L^{\Phi}\left(G^{u}\right)$ for every $u \in G^{0}$ and using lemma $3.4$ in the context of $M^{\Psi}$, if $S \in R\left(E^{\Psi}, C_{0}(G^{0})\right)$, there exist a unique $g \in E_{0}^{\Phi}$ such that $S(\xi)(u)=\langle g, \xi\rangle(u)$. Also,
\begin{align*}
\left\|g^{u}\right\|_{\Phi}^{0} & =\sup \left\{|\int g^{u}fd\lambda^{u}|: f \in M^{\Psi}\left(G^{u}\right),\|f\|_{\Psi} \leq 1\right\}, \\
& =\sup \left\{|S(\xi)(u)|: \xi \in E^{\Psi},\|\xi\|_{\Psi} \leq 1\right\} .
\end{align*}
So $\|g\|_{\Phi}^{0}=\|S\|$. Hence we can see that $R\left(E^{\Psi}, C_{0}(G^{0})\right)$ is closed in $B_{D}\left(E^{\Psi}, C_{0}(G^{0})\right)$.

Let $\check{A}_{\Phi}(G)^{\prime}$ be the set of $\alpha \in B_{D}^{\Phi}\left(\check{A}_{\Phi}(G), C_{0}(G^{0})\right)$, such that $f \alpha \in  R\left(E^{\Psi}, C_{0}(G^{0})\right)$ for all $f \in C_{c}(G)$. It can be easily verified that $\check{A}_{\Phi}(G)^{\prime}$ is a closed subspace of $B_{D}^{\Phi}\left(\check{A}_{\Phi}(G), C_{0}(G^{0})\right)$.

The following is the main theorem of this section. For the proof we use some ideas from \cite{paterson2004fourier}
and \cite{aghababa2013space}.
\begin{thm}
     Let $G$ be a groupoid, $(\Phi, \Psi)$ be a complementary pair of $N$-functions both satisfying $\Delta_{2}$-condition. Then $\check{A}_{\Phi}(G)^{\prime}$ can be identified with $C_{V_{\Phi}}(G)$.
\end{thm}
\begin{proof}
    Let $T \in C_{V_{\Phi}}(G), h \in \check{A}_{\Phi}(G)$. Then for some $P \in \mathcal{K}(G)$,$$ h=\sum_{i=1}^{\infty} g_{i} * \check{f}_{i},\ f_{i}\in E^{\Phi}_{0}(P), g_{i} \in E^{\Psi}_{0}(P).$$ 
    Define $\phi_{T}: \check{A}_{\Phi}(G) \rightarrow C_{0}\left(G^{0}\right)$ by
\[\phi_{T}(h)(u)=\sum_{i=1}^{\infty}\left\langle T f_{i}, g_{i}\right\rangle(u).\]
It is clear that $\phi_{T}$ is linear. Further, $$\left\|\phi_{T}(h)\right\|_{\infty} \leq 2 \sum_{i=1}^{\infty}\left\|T f_{i}\right\|_{\Phi}^{0}\left\|g_{i}\right\|_{\Psi}^{0} \leq 2\|T\| \sum_{i=1}^{\infty}\left\|f_{i}\right\|_{\Phi}^{0}\left\|g_{i}\right\|_{\Psi}^{0}<\infty.$$ We need to show that $\phi_{T}(h)$ is independent of the representation of $h$. For that, it is enough to show that, $\phi_{T}(h)=0$, if $h=0$. Suppose $h=\sum_{i=1}^{\infty} g_{i} * \check{f}_{i}=0$. Also, $T_{k}(f) =T\left(e_{k} * f\right)= Te_{k}*f, ~ f\in C_{c}(G) $. For all $u \in G^{0}$,
\begin{align*}
\left|\sum_{i=1}^{\infty}\left\langle T_{k} f_{i}, g_{i}\right\rangle(u)\right|  \leq 2 \sum_{i=1}^{\infty}\left\|T_{k} f_{i}\right\|_{\Phi}^{0}\left\|g_{i}\right\|_{\Psi}^{0}&
 \leq 2\left\|T_{k}\right\| \sum_{i=1}^{\infty}\left\|f_{i}\right\|_{\Phi}^{0}\left\|g_{i}\right\|_{\Psi}^{0} \\
& \leq 4\|T\| \sum_{i=1}^{\infty}\left\|f_{i}\right\|_{\Phi}^{0}\left\|g_{i}\right\|_{\Psi}^{0}<\infty.
\end{align*}
So, using Lemma $4.1$, the function $\sum_{i=1}^{\infty}\left\langle T_{k} f_{i}, g_{i}\right\rangle$ converges uniformly in $k$ to $\sum_{i=1}^{\infty}\left\langle Tf_{i}, g_{i}\right\rangle$ on $C_{0}\left(G^{0}\right)$.

Let $n_{j} \in \mathbb{N}$ such that  $\sum_{n>n_{j}}\|f_n\|^{0}_{\Phi}\|g_{n}\|_{\Psi}^{0} < \frac{1}{2jM}$ where $M = \max \{4\|T\|, 2\}$. Since $\{f\in C_{c}(G): \operatorname{supp}(f)\subset P\}$ is dense in $E^{\Phi}_{0}(P)$, there exist a family of sequences $\left\{(h_{i}^{j})_{i\in \mathbb{N}}: h^{j}_{i}\in C_{c}(G), \operatorname{supp}(h^{j}_{i}) \subset P,~j\in \mathbb{N}\right\}$  such that $h^{j}_{i}=0 ~\forall~ i>n_{j}$ and $\|f_{i}-h^{j}_{i}\|^{0}_{\Phi}< \frac{1}{2jK_{j}M}$, for $1\leq i\leq n_{j}$ where $K_{j}= \sum_{i=1}^{n_{j}}\|g_{i}\|^{0}_{\Psi}$.
Hence,

\begin{align*}
\sum_{i=1}^{\infty}\left\langle T_{k} f_{i}, g_{i}\right\rangle(u)  & =\lim_{j\to \infty}\sum_{i=1}^{n_{j}}\left\langle T e_{k} * h_{i}^{j}, g_{i}\right\rangle(u)
\\ & =\lim_{j\to \infty}\sum_{i=1}^{n_{j}}\left\langle\chi_{P^*} T e_{k}, g_{i} * \check{h^{j}_{i}}\right\rangle(u) \\ & =\left\langle\chi_{P^*} T e_{k}, \lim_{j\to \infty}\sum_{i=1}^{n_{j}} g_{i} * \check{h_{i}^{j}}\right\rangle(u)\\
&= \left\langle\chi_{P^*} T e_{k}, \sum_{i=1}^{\infty} g_{i} * \check{f_{i}}\right\rangle(u)
\end{align*}
Here $P^{*}=P P^{-1}$ is compact in $G$. Note that $ \chi_{P^{*}} T e_{k}$ denotes the section whose value at $u$ is the function: $$\left(\chi_{P^{*}} T e_{k}\right)^{u}(t)=\chi_{P^{*} \cap G^{u}}(t)\left(T e_{k}\right)^{u}(t)$$ for $t \in G^{u}.$ We can easily see that $\chi_{P^{*}} T e_{k}\in I(G,\lambda)$. So, $\sum_{n=1}^{\infty}\left\langle T f_{n}, g_{n}\right\rangle(u)=0$, for all $u \in G^{0}$. Hence $\phi_{T}(h)$ is well defined and $\left\|\phi_{T}(h)\right\|_{\infty} \leq 2\|T\|\|h\|_{\check{A}_{\Phi(G)}}.$ Further, 
\begin{align*}
\|T\| & =\sup \left\{\left\|(T f)^{u}\right\|_{\Phi}^{0}: u \in G^{0}, f \in C_{c}(G),\|f\|_{\Phi}^{0} \leq 1\right\}, \\
& =\sup \left\{|\langle T f, g\rangle(u)|: u \in G^{0},  g ,f \in C_{c}(G),\|f\|_{\Phi}^{0} \leq 1,\|g\|_{\Psi} \leq 1\right\}, \\
& =\sup \left\{\left\|\phi_{T}(h)\right\|_{\infty}: h=g * \check{f},\|h\|_{\check{A}_{\Phi(G)}}~ \leq 1\right\}
\leq\left\|\phi_{T}\right\| .
\end{align*}
So $\|T\| \leq\left\|\phi_{T}\right\| \leq 2\|T\|$. For $b \in C_{0}\left(G^{0}\right), h \in \check{A}_{\Phi(G)}$,
\begin{align*}
\phi_{T}(h b)(u) =\sum_{n=1}^{\infty}\left\langle T f_{n}, g_{n} b\right\rangle(u)&=\sum_{n=1}^{\infty}\left\langle T f_{n}, g_{n}\right\rangle(u) b(u) 
 =\left(\phi_{T}(h) b\right)(u) .
\end{align*}
Thus $\phi_{T}$ is a bounded linear right $C_{0}(G^{0})$-module map. Now for $f \in C_{c}(G),$ $$f \phi_{T}(g)(u)=\phi_{T}(g * \check{f})(u)=\langle T f, g\rangle(u),$$ for all $g \in C_{c}(G).$ Hence $f \phi_{T} \in R\left(E^{\Psi}, C_{0}(G^{0})\right)$ implies $\phi_{T} \in \check{A}_{\Phi}(G)^{\prime}.$ 

It remains to show that $\phi: T \rightarrow \phi_{T}$ is surjective. Let $\alpha \in \check{A}_{\Phi}(G)^{\prime}$, then for $f \in C_{c}(G)$, there exist $F_{f} \in E_{0}^{\Phi}$ such that, $$(f \alpha)(g)(u)=\left\langle F_{f}, g\right\rangle(u), g \in E^{\Psi}.$$ Define $T(f)=F_{f}$, for all $f \in C_{c}(G).$ Note $\left\langle F_{f}, g\right\rangle(u)=\langle T(f), g\rangle(u)=\alpha(g * \check{f})(u)$, for all $g \in C_{c}(G).$ Further, $$|\langle T(f), g\rangle(u)|=|\alpha(g * \check{f})(u)| \leq\|\alpha\|\|f\|_{\Phi}^{0}\|g\|_{\Psi}.$$ Hence $g \in C_{c}(G)$ implies $\|T f\|_{\Phi}^{0} \leq\|\alpha\|\|f\|_{\Phi}^{0}$. So, $T$ can be extended to $E_{0}^{\Phi}$ and $\|T\| \leq\|\alpha\|$. Let $f_{1}, f_{2}, g \in C_{c}(G)$. Then, 
\begin{align*}
\langle T(f_{1} * f_{2}), g\rangle & = \alpha(g *\check{(f_{1} * {f}_{2}})) \\ & = \alpha(g *(\check{f}_{2} * \check{f}_{1})) \\
& =\alpha\left(\left(g * \check{f}_{2}\right) * \check{f}_{1}\right) \\
& =\left\langle T f_{1}, g * \check{f}_{2}\right\rangle =\left\langle T f_{1} * f_{2}, g \right 
\rangle .
\end{align*}
 Hence, $T \in C_{V_{\Phi}}(G)$.
\end{proof}

\section*{Acknowledgement}
 K. N. Sridharan is supported by  NBHM doctoral fellowship with Ref number: 0203/13(45)/2021-R\&D-II/13173. The authors appreciate the editor and the anonymous reviewer for their careful study of the manuscript and for offering valuable suggestions that enhanced the quality of the article.

\section*{Data Availability}
Data sharing does not apply to this article as no datasets were generated or analyzed during the current study.

\section*{competing interests}
The authors declare that they have no competing interests.

\small{\bibliographystyle{unsrt}}
\bibliography{reference}

\end{document}